\documentclass[12pt,twoside,a4paper]{amsart}
\usepackage{amssymb}
\date{\today}

\begin{document}

\title[]
{A generalized Poincar\'e-Lelong formula}

\def\la{\langle}
\def\ra{\rangle}
\def\End{{\rm End}}

\def\deg{\text{deg}\,}


\def\al{\text{\tiny $\aleph$}}
\def\w{\wedge}
\def\db{\bar\partial}
\def\dbar{\bar\partial}
\def\ddbar{\partial\dbar}
\def\loc{{\rm loc}}
\def\R{{\mathbb R}}
\def\C{{\mathbb C}}
\def\w{{\wedge}}
\def\P{{\mathbb P}}
\def\cp{{\mathcal C_p}}
\def\bl{{\mathcal B}}
\def\A{{\mathcal A}}
\def\B{{\mathcal B}}
\def\Cn{\C^n}
\def\ddc{dd^c}
\def\bbox{\bar\square}
\def\G{{\mathcal G}}
\def\od{\overline{D}}
\def\ot{\leftarrow}
\def\loc{\text{loc}}
\def\D{{\mathcal D}}
\def\M{{\mathcal M}}
\def\Hom{{\rm Hom\, }}
\def\codim{{\rm codim\,}}
\def\Tot{{\rm Tot\, }}
\def\Im{{\rm Im\, }}
\def\K{{\mathcal K}}
\def\Ker{{\rm Ker\,  }}
\def\Dom{{\rm Dom\,  }}
\def\Z{{\mathcal Z}}
\def\E{{\mathcal E}}
\def\can{{\rm can}}
\def\O{{\mathcal O}}
\def\Pop{{\mathcal P}}
\def\L{{\mathcal L}}
\def\Q{{\mathcal Q}}
\def\Re{{\rm Re\,  }}
\def\Res{{\rm Res\,  }}
\def\Aut{{\rm Aut}}
\def\L{{\mathcal L}}
\def\U{{\mathcal U}}

\def\sr{\stackrel}

\newtheorem{thm}{Theorem}[section]
\newtheorem{lma}[thm]{Lemma}
\newtheorem{cor}[thm]{Corollary}
\newtheorem{prop}[thm]{Proposition}

\theoremstyle{definition}

\newtheorem{df}{Definition}

\theoremstyle{remark}

\newtheorem{preremark}{Remark}
\newtheorem{preex}{Example}

\newenvironment{remark}{\begin{preremark}}{\qed\end{preremark}}
\newenvironment{ex}{\begin{preex}}{\qed\end{preex}}

\numberwithin{equation}{section}

\date{\today}

\author{Mats Andersson}

\address{Department of Mathematics\\Chalmers University of Technology and the University of G\"oteborg\\
S-412 96 G\"OTEBORG\\SWEDEN}

\email{matsa@math.chalmers.se}



\thanks{The author was
  partially supported by the Swedish Natural Science
  Research Council}

\begin{abstract} We prove a  generalization of the
classical  Poincar\'e-Lelong formula. 
Given a holomorphic section $f$, with zero set $Z$, of  a Hermitian 
vector bundle $E\to X$, let
$S$ be the line bundle over $X\setminus Z$ spanned by $f$ and let $Q=E/S$. Then 
the Chern form $c(D_Q)$ is locally integrable and closed in $X$ and  there is a current
$W$ such that
$dd^cW=c(D_E)-c(D_Q)-M,$
where $M$  is a current with support on $Z$. 
In particular, 
the top Bott-Chern class is represented by a current with support on $Z$.
We discuss positivity of these  currents,  and
we also reveal a close  relation with  principal value 
and residue currents of Cauchy-Fantappi\`e-Leray type.

\end{abstract}

\maketitle

\section{Introduction}

Let $f$ be a holomorphic (or meromorphic) section of  a  Hermitian  line bundle
$L\to X$, and let $[Z]$ be the current of integration over the 
divisor $Z$ defined by $f$.  The Poincar\'e-Lelong formula  states that
$$
dd^c\log(1/|f|)=c_1(D_L)-[Z],
$$
where $c_1(D_L)$ is the first Chern form associated with the Chern connection
$D_L$ on $L$, i.e.,
$c_1(D_L)=\al\Theta_L$, where $\Theta_L$ is the curvature;
here and throughout this paper $\al=i/2\pi$ 
and 
$d^c=\al(\dbar-\partial)$ so that
$$
dd^c=\frac{i}{\pi}\partial\dbar=2\al\partial\dbar.
$$
If  $U$ is  the meromorphic section of the dual bundle $L^*$ such that
$U\cdot f=1$, then $R=\dbar U$ is a $(0,1)$-current, and we have the global
factorization 
\begin{equation}\label{fakt}
[Z]= R\cdot D_Lf/2\pi i.
\end{equation}
If
$A=-2\al\partial \log(1/|f|)$, then clearly
$dA=\dbar A=c_1(D_L)-[Z]$, and it is easily checked that
$A= U\cdot D_Lf/2\pi i$.  
In this paper we consider analogous formulas
for a holomorphic section $f$ of  a higher rank bundle, and our main result
is the following generalization of
the  Poincar\'e-Lelong formula.

\begin{thm}\label{thmett}
Let $f$ be a holomorphic section of   the Hermitian  vector bundle
$E\to X$ of rank $m$. Let $Z=\{f=0\}$, let $S$ denote the (trivial) line bundle over
$X\setminus Z$ generated by $f$, and let $Q=E/S$, equipped with the induced 
Hermitian metric. 
\smallskip

\noindent(i)
The Chern form $c(D_Q)$ is locally integrable in $X$ and its natural
extension to $X$ is  closed.    Moreover, the 
 forms  $\log|f|c(D_Q)$ and 
\begin{equation}\label{elefant}
|f|^{2\lambda}\frac{\al\partial|f|^2\w\dbar|f|^2}{|f|^4}\w c(D_Q),\quad \lambda>0,
\end{equation}
are locally integrable in $X$,  and
\begin{equation}\label{Md}
M=\lim_{\lambda\to 0^+}\lambda|f|^{2\lambda}\frac{\al\partial|f|^2\w\dbar|f|^2}{|f|^4}\w c(D_Q)
=dd^c(\log|f| c(D_Q)){\bf 1}_Z
\end{equation}
is a closed current of order zero with support on $Z$. If $\codim Z=p$, then
$$
M=M_p+M_{p+1}\cdots +M_{\min(m,n)},
$$
where $M_k$ has bidegree $(k,k)$, and
$$
M_p=\sum\alpha_j[Z_j^p],
$$
where $Z^p_j$ are the irreducible components of codimension  precisely $p$,  and
$\alpha_j$ are the Hilbert-Samuel multiplicities of  $f$.

\smallskip
\noindent(ii)
There is a current $W$ of bidegree $(*,*)$ and order zero in $X$  which is  smooth
in $X\setminus Z$, 
and with logarithmic singularity at $Z$,  such that
\begin{equation}\label{bas}
dd^cW=c(D_E)-C(D_Q)-M,
\end{equation}
where $C(D_Q)$ denote the natural extension of $c(D_Q)$.
\end{thm}

Here $c(D)$ denotes the Chern form with respect to the Chern connection $D$ 
associated to the Hermitian  structure, i.e.,
$c(D)=\det(\al\Theta+I)$, 
where $\Theta=D^2$ is the curvature tensor.
We let $c_k(D)$ denote the component of bidegree $(k,k)$.

\smallskip
For an  explicit expression for $W$, see Definition~\ref{plask} in
Section~\ref{grotta}.  
If  $W_k$ denotes  the  component of bidegree $(k,k)$,  then
\eqref{bas} means that
\begin{equation}\label{bas2}
dd^cW_{k-1}=c_k(D_E)-c_k(D_Q)-M_k.
\end{equation}
Since $Q$ has rank $m-1$,  $c_m(D_Q)=0$, and therefore 
$$
dd^cW_{m-1}=c_m(D_E)-M_m,
$$
which means that the current $M_m$ represents the top degree Bott-Chern
 class $\hat c_m(E)$.
It also follows that the Bott-Chern class $\hat c_k(E)$ is equal to
$\hat c_k(Q)$ if $k<p$.

\smallskip
If  $E$ is a line bundle, then, see Definition~\ref{plask},
$W=W_0=\log(1/|f|)$,  so \eqref{bas2} is the
then  usual Poincar\'e-Lelong formula.

\smallskip
In \cite{BC} Bott and Chern developed a method of transgression
which  in particular gives a form $w$ in $X\setminus Z$ such that
$dd^c w=c(D_E)-c(D_Q)$.
It is not  unexpected that one can extend this construction
across $Z$ by a careful analysis  of the occurring singularities
at $Z$.
In the recent paper \cite{Meo2}, Meo proves \eqref{bas2} for
$k=p$.  Previously this formula was proved
in \cite{BGS} in the case when $f$ defines a complete intersection, i.e.,
$p=m$.  A variety of analogous formulas
for $d$ rather than $dd^c$ are constructed in 
quite  general (non-holomorphic)  situations  in \cite{HS},
\cite{HL1}, \cite{HL2}, and \cite{HL3}.
\smallskip

Clearly $M_p$ is always a positive current. It  follows from \eqref{Md}  that
$M_k$ is positive if $c_{k-1}(D_Q)$ is a positive form. 
For an even more precise formula for $M$, see
Proposition~\ref{mprecis}.

Let us say that $E$ is positive if $E^*$ is Nakano negative.

\begin{thm}\label{thmettx}
Assume that $E$ is positive. Then $c(D_E)$ is a positive form, 
$C(D_Q)$ and $M$ are positive currents,
and  (one can choose $W$ such that)
$W$ is positive where $|f|\le 1$.
\end{thm}

If $A=-2\al\partial W$ we have, cf., \eqref{bas},
\begin{equation}\label{Apa}
\dbar A=dA=c(D_E)-c(D_Q)-M.
\end{equation}
In \cite{A2} we introduced  a residue current
 $R=R_p+\cdots +R_{\min(m,n)}$,
associated with $f$,
with support on $Z$, where $R_k$ is a $(0,k)$-current with
values in $\Lambda^k E^*$,
and a principal value current  $U=U_1+\cdots +U_m$
such that 
$(\delta_f-\dbar)U=1-R$, 
where $\delta_f$ denotes contraction with $f$.
When $E$ is a line bundle, then
$U=1/f$ and $R=\dbar(1/f)$.
In analogy to \eqref{fakt} we can factorize $M_p$ as
$$
M_p=R_p\cdot (D_Ef)^p/p!;
$$
this was proved  in \cite{A3}.
We have a similar, but somewhat more involved, formula for
the whole current $M$, see \eqref{gris}  in Section~\ref{faktor}.
In a similar  way we can express $A$ and $c(D_Q)$, see \eqref{gris2} and \eqref{cformel}, 
in terms of the current $U$.

\begin{remark}
Let  $f_1,\ldots,f_r$ be holomorphic sections of $E$ and let $Z$ be the 
analytic set where they are linearly dependent. 
Moreover, let $S$ be the trivial rank $r$-subbundle of $E$ over $X\setminus Z$
generated by $f_j$ and let $Q=E/S$. Then $c(D_Q)$ has a natural current extension $C(D_Q)$
across $Z$  and there is a closed current $M$ of bidegree  $(*,*)$ with support
on $Z$ and a current $A$  such that
\begin{equation}\label{anita}
dA=\dbar A=c(D_E)-C(D_Q)-M.
\end{equation}
This can be proved by a small modification of the argument in this paper;
in   the case $Z$ has generic dimension such a formula was  proved already in
\cite{HL2}, and the general case should be  contained in \cite{HL3}.
It follows from \eqref{anita} that the   current
$M_k$ is a representative for $c_k(D_E)$ for $k>m-r$.

However,  we have no analogous formula  for $dd^c$.
\end{remark}

As indicated above, the proof of Theorem~\ref{thmett}
relies on the construction in \cite{BC},
combined with a careful control of the singularities
at $Z$.
To begin with one constructs  a form $v$  in $X\setminus Z$  such that
$$
dd^c v=c(D_E)-c(D_S)c(D_Q).
$$
By Hironaka's theorem and toric resolutions, following \cite{BGVY} and \cite{PTY}, 
we  can prove that this  equality has  meaning in the current sense across
$Z$.  Here a  crucial point is  an explicit formula for the Chern form $c(D_Q)$
(Proposition~\ref{cqprop}) from  which it is easy to conclude that 
$c(D_Q)$  has a smooth extension across the singularity after an appropriate blow-up.
By the  usual Poincar\'e-Lelong formula, 
$c(D_S)-1=dd^c\log(1/|f|)$ outside $Z$, and we can conclude that
\eqref{bas} holds  (if the capitals denote  the natural extensions across $Z$) with
$$
W=   \log(1/|f|)C(D_Q)-V,
$$
and
$M=dd^c(\log|f|C(D_Q)){\bf 1}_Z.$
Theorem~\ref{thmettx} follows  essentially  by applying ideas in \cite{BC}.

In Section~\ref{poss} we discuss the positivity
and prove Theorem~\ref{thmettx}, essentially by applying ideas from \cite{BC}.
The paper is concluded by some  examples.

\section{Preliminaries}

We first  recall
the differential geometric  definition of  Chern classes.
Let $E\to X$ be any differentiable complex vector bundle over
a differential manifold $X$, with connection
$D\colon\E_k(X,E)\to\E_{k+1}(X,E)$ 
and curvature tensor
$D^2=\Theta\in\E_2(X,\End  E)$.
The connection $D=D_E$ induces in a natural way a connection
$D_{\End  E}$ on the bundle $\End E$ by the formula 
$Dg\cdot\xi=D(g\cdot \xi)-g\cdot D\xi$,
and in a similar way there is a natural  connection $D_{E^*}$ on the dual bundle  
$E^*$, etc.
In particular we have Bianchi's identity
\begin{equation}\label{bianchi}
D_{\End E}\Theta=0.
\end{equation}
If  $I$ denotes the identity mapping on $E$, then
$c(D)=\det(\al\Theta +I)$
is a welldefined differential form whose terms have  even degrees, which is
 called the Chern form of $D$. It is a  basic fact that
$c(D)$  is a closed form.
Moreover  its de~Rham cohomology class
 is independent of $D$ and is called the (total) 
Chern class $c(F)$  of the bundle $F$.

\smallskip

To prove this, one can  consider a smooth one-parameter family $D_t$
of connections of $F$ with $D_0=D$. If $E'$ is the pull-back of $E$
to  $X\times[0,1]$, then 
$D'=D_t+d_t$
is a connection on $E'$ and its curvature tensor is
$$
\Theta'=\Theta_t+dt\wedge \dot{D}_t
$$
where
$\dot{D}_t=d D_t/dt$. It is readily checked that it is 
an element in $\E_1(X,\End(F))$. 
Since
$(d+d_t)\det(\al\Theta'+I)=0$
we have that
$$
d_\zeta\int_0^1\det(\al\Theta'+I)=
-\int_0^1 d_t\det(\al \Theta'+I)=
c(D)-c(D_1).
$$

\smallskip
In order to make the  computation more explicit we
introduce the exterior algebra bundle
$\Lambda=\Lambda(T^*(X)\oplus F\oplus F^*)$. Any 
section $\xi\in\E_k(X,F)$ corresponds to a section $\tilde\xi$ of $\Lambda$;
if $\xi=\xi_1\otimes e_1+\ldots +\xi_m\otimes e_m$ in a local frame
$e_j$, then we let 
$\tilde\xi=\xi_1 \w e_1+\ldots +\xi_m \w e_m$.
In the same way,   $a\in\E_k(X,\End E)$ can be identified with
$$
\widetilde a=\sum_{jk} a_{jk} \wedge e_j\wedge e^*_k,
$$
if $e_j^*$ is the dual frame,
 and $a=\sum_{jk} a_{jk}\otimes e_j\otimes e^*_k$ 
with respect to these frames.
A given connection  $D=D_F$ on $F$ extends 
in a unique  way to a linear  mapping
$\E(X,\Lambda)\to\E(X,\Lambda)$ which is
a an anti-derivation with respect to the
wedge product in $\Lambda$, and such that it acts as
the exterior differential $d$ on the $T^*(X)$-factor.
It is readily seen that
$$
\widetilde{D_E\xi}=D\tilde\xi,
$$
if $\xi$ is a form-valued section of  $E$. In the same way we have
\begin{lma}\label{lapp}
If $a\in\E_k(X,\End E)$, then
\begin{equation}
\widetilde{D_{\End E}a}=D\widetilde a.
\end{equation}
\end{lma}

\begin{proof}
If $\xi\in\E_k(X,E)$ and $\eta\in\E(X,E^*)$, then
$$
D_{\End E}(\xi\otimes\eta)=D_E\xi\otimes\eta+(-1)^k\xi\otimes D_{E^*}\eta,
$$
and thus the snake of $D_{\End E}(\xi\otimes\eta)$ is equal to
$$
\widetilde{D_E\xi}\w\eta+(-1)^{k+1}\tilde\xi\wedge\widetilde{ D_{E^*}\eta}
=D(\tilde\xi\wedge\eta)
$$
as claimed.
\end{proof}

Since $D_{\End E}I=0$,   ($I=I_E$) we have from \eqref{bianchi}
and Lemma~\ref{lapp} that 
\begin{equation}\label{enol}
D\widetilde\Theta=0 \quad \text{and}\quad  D\tilde I=0.
\end{equation}
We let  $\tilde I_m=\tilde I^m/m!$  and  use the same notation 
for other forms in the sequel.
Any form $\omega$ with values in $\Lambda$ can be  written 
$\omega=\omega'\w\tilde I_m+\omega''$ uniquely, 
where $\omega''$ has lower  degree in $e_j,e^*_k$. If we define
$$
\int_e\omega =\omega',
$$
then this integral is of course linear and moreover
\begin{equation}\label{skodon}
d\int_e\omega=\int_e D\omega.
\end{equation}
In fact, since $D\tilde I_m=0$, 
$$
\int_e D\omega=\int_e d\omega'\wedge \tilde I_m+D\omega''=d\omega'=d\int_e\omega.
$$
Observe that
\begin{equation}\label{laban}
c(D)=\int_e(\al\widetilde\Theta+\tilde I)_m=
\int_e e^{\al\tilde\Theta+\tilde I}.
\end{equation}
Lemma~\ref{lapp} and  \eqref{enol} together imply that
the Chern form $c(D)$ is closed.  Furthermore, following the outline
above,  we get the formula
\begin{equation}\label{smuts}
d\int_0^1\int_e\al
\widetilde{\dot{D}}\wedge
e^{\al\widetilde\Theta_t+\tilde I}
=c(D_1)-c(D_0),
\end{equation}
thus showing that $c(D_0)$ and $c(D_1)$ are cohomologous.

\smallskip
Recall that if the connection $D$ is modified to $D_1=D-\gamma$, 
where $\gamma\in\E_1(X,\End E))$,  then
$\Theta_1=\Theta -D_{\End E}\gamma +\gamma\w\gamma$.
If we form the explicit  homotopy
$D_t=D- t\gamma$, therefore
\begin{equation}\label{vatten}
\Theta_t=\Theta-tD_{\End E}\gamma+t^2\gamma\wedge\gamma
\end{equation}
and hence, by Lemma~\ref{lapp},
\begin{equation}\label{bror}
\widetilde\Theta_t=\widetilde\Theta-tD\tilde\gamma+t^2\widetilde{\gamma\w\gamma}.
\end{equation}

\section{Bott-Chern classes}\label{bcsektion}

From now on we assume that $E$ is a holomorphic
Hermitian  bundle and that $D_E$ is the Chern connection
and $D'_E$ is its  $(1,0)$-part. Then the induced connection
$D_{E^*}$ on $E^*$ is the Chern connection on $E^*$ etc.
In particular, our mapping $D$ on $\Lambda$ is of type
$(1,0)$, i.e., $D=D'+\dbar$.

\smallskip

Let  $E\to X$ be   a Hermitian  vector bundle with Chern connection $D_E$. The
Bott-Chern class $\hat c(E)$ is the equivalence class of the Chern form $c(D_E)$ in
$$
\frac{\oplus_k\E_{k,k}(X)\cap \Ker d}{\oplus_kdd^c\E_{k,k}(X)}.
$$

\smallskip
\begin{lma} Let $D$ be a connection depending smoothly on a real parameter $t$.
Moreover, assume that $L\in\E(X,\End(E))$ depends smoothly on $t$ and that
\begin{equation}\label{cbvillkor}
D'_{\End E}L= \dot{D}.
\end{equation}
Also  assume that $\Theta_t$ has bidegree
$(1,1)$ for all $t$.
If 
$$
v=-\frac{1}{2}\int_0^1\int_e\tilde L_t\w e^{\al\tilde\Theta_t+\tilde I}dt,
$$
then $-2\al\partial v=b$, where
$$
b=\int_0^1\int_e\al\widetilde{\dot{D}_t} \w e^{\al\tilde\Theta_t+\tilde I}dt.
$$
\end{lma}

This lemma as well as the other material in this section is
taken from \cite{BC}. However, we use a somewhat different
formalism, and for the reader's convenience we supply
some simple proofs.

\begin{proof}
In view of \eqref{skodon} 
we have that (suppressing the index $t$)
$$
d\int_e\tilde L\w e^{\al\tilde\Theta+\tilde I}=
\int_e D\tilde L\w e^{\al\tilde\Theta+\tilde I},
$$
and by identifying bidegrees we get that
$$
\partial \int_e\tilde L\w e^{\al\tilde\Theta+\tilde I}=
\int_e D'\tilde L\w e^{\al\tilde\Theta+\tilde I}=
\int_e \widetilde{\dot{D}} \w e^{\al\tilde\Theta+\tilde I}.
$$
\end{proof}

Since $db=c(D_1)-c(D_0)$, cf.,  \eqref{smuts}, 
we thus have
\begin{equation}\label{cbhom}
-dd^c v=c(D_1)-c(D_0).
\end{equation}

\smallskip
By deforming the metric one can  use this lemma to show that
$\hat c(E)$ is independent of the Hermitian  structure on $E$, see \cite{BC}.
However we are interested in a somewhat different   situation. 
Assume that we  have  the short exact sequence of Hermitian  vector bundles
\begin{equation}\label{short}
0\to S\sr{j}{\rightarrow} E\sr{g}{\rightarrow} Q\to 0,
\end{equation}
where $Q$ and $S$ are equipped with the metrics induced by the Hermitian  metric of $E$.
Then 
\begin{equation}\label{iso} 
j^*\oplus g \colon E\to S\oplus Q
\end{equation}
is a smooth vector bundle isomorphism. If $D_S$ and $D_Q$ are the Chern connections
on $S$ and $Q$ respectively, then
\begin{equation}\label{betaox}
D_E\sim
\left(\begin{array}{cc}
D_S & -\beta^* \\
\beta & D_Q 
\end{array}
\right)
\end{equation}
with respect to the isomorphism \eqref{iso},
where $\beta\in\E_{1,0}(X,\Hom(S,Q))$  is the second fundamental form, see \cite{Dem}. 
We shall now modify the connection $D=D_E$ to
$D_b=D-\gamma_b$, where
$\gamma_b=D'_{\End E}jj^*$.
It turns out that  $\gamma=g^*\circ\beta\circ j^*$, thus 
$\gamma\w\gamma=0$, and that 
$D_{\End E}\gamma=\dbar\gamma$.
Moreover, it follows that
$$
D_{b}\sim
\left(\begin{array}{cc}
D_S &  * \\
0  & D_Q 
\end{array}
\right)
$$
and hence
\begin{equation}\label{thb}
\Theta_b\sim
\left(\begin{array}{cc}
\Theta_S &  * \\
0  & \Theta_Q 
\end{array}
\right),
\end{equation}
so that $c(D_b)=c(D_S)c(D_Q)$.
If $D_{t}=D-t\gamma_b$ we have 
$\Theta_t=\Theta-t\dbar\gamma_b$;  thus it has bidegree $(1,1)$.
If we let
\begin{equation}\label{gaas}
b=\int_0^1\int_e \al\tilde\gamma_b\wedge e^{\tilde I+\al\tilde\Theta-t\al\dbar\tilde\gamma_b}
=\sum_{\ell\ge 0}\int_e\al\tilde\gamma_b\w e^{\tilde I+\al\widetilde\Theta}\w\frac{1}{(\ell+1)!}
(-\al\dbar\tilde\gamma_b)^{\ell}
\end{equation}
it follows from \eqref{smuts} that
$d b=c(D_E)-c(D_S)c(D_Q)$.
Moreover, if $L=jj^*/(1-t)$, 
then \eqref{cbvillkor} holds. In fact, $\dot D=-\gamma_b$, 
and $[jj^*,g^*\circ\beta\circ j^*]=g^*\circ\beta\circ j^*,$
so that
\begin{equation}\label{ramaskri}
D'_{\End E,t}L=D'_{\End E}L-t[\gamma_b,L]=
\frac{1}{1-t}\gamma_b-\frac{t}{1-t}\gamma_b=\gamma_b.
\end{equation}


\begin{prop}\label{bc}
If 
\begin{equation}\label{bla}
v=\sum_{\ell=1}^{m-1}\frac{(-1)^{\ell}}{2\ell}\int_e\widetilde{jj^*}\wedge
(\tilde I+\al\tilde\Theta-\al\dbar\tilde\gamma_b)_{m-\ell-1}\wedge
(-\al\dbar\tilde\gamma_b)_\ell,
\end{equation}
then
$-2\al\partial v=b.$
\end{prop}


\begin{proof}
Observe that
$$
\partial\int_0^{1-\epsilon}\int_e\frac{\widetilde{jj^*}}{1-t}\w 
e^{\tilde I+\al\Theta_1}dt=
\int_0^{1-\epsilon}\int_e\frac{D_1\widetilde{jj^*}}{1-t}\w 
e^{\tilde I+\al\Theta_1}dt=0,
$$
since $D_1\widetilde{jj^*}=\widetilde{D_{\End E,1}jj^*}=0$ in view of Lemma~\ref{lapp} and \eqref{ramaskri}. Therefore,
$$
\al\partial\int_0^{1-\epsilon}\int_e \widetilde{jj^*}\wedge\frac{e^{\tilde I+\al\tilde\Theta
-t\al\dbar\tilde\gamma_b}-e^{\tilde I+\al\tilde\Theta
-\al\dbar\tilde\gamma_b}}
{1-t}dt=\int_0^{1-\epsilon}
\int_e \al\tilde\gamma_b\wedge e^{\tilde I+\al\tilde\Theta-t\al\dbar\tilde\gamma_b}.
$$
The proposition now follows by letting  $\epsilon\to 0$ and computing
the $t$-integral on the left hand side.
\end{proof}
Altogether we therefore have that
$-dd^c v=c(D_E)-c(D_S) c(D_Q)$ and thus
$\hat c(E)=\hat c(S)\hat c(Q)$. 

\section{Proof of the main formula}\label{grotta}

Let $f$ be a nontrivial  holomorphic section of $E$,  $Z=\{f=0\}$, and let
$S$ be  the  trivial subbundle of $E$ over $X\setminus Z$, 
generated by the $f$.
We then have the  short exact sequence \eqref{short}
over $X\setminus Z$, where $g\colon E\to Q=E/Q$ is the natural projection.
Let $\sigma$ be the   section  of  the dual bundle $E^*$
with minimal norm such   $\sigma \cdot f=1$. Then clearly
\begin{equation}\label{projj}
\widetilde{jj^*}= f\w\sigma.
\end{equation}

Observe that the natural  conjugate-linear isometry
$E\simeq E^*$, $\eta\mapsto \eta^*$, defined by
$$
\eta^*\cdot\xi=\langle\xi,\eta\rangle,\quad \xi\in\E(X,E), 
$$
extends to an isometry on the space of form-valued sections.

\begin{lma}\label{lmalma} 
If $\phi=-\partial\log|f|^2$, then 
$D'\sigma=\phi\wedge\sigma.$
\end{lma}

\begin{proof}
Observe that $\sigma=f^*/|f|^2$. Since $D=D_E$ is the Chern connection,
 $D' f^*=(\dbar f)^*=0$,   so we have 
$$
D'\sigma=D'(f^*/|f|^2)=\partial\frac{1}{|f|^2}\wedge f^*=-\partial\log|f|^2\wedge\sigma.
$$
\end{proof}

Following Section~\ref{bcsektion}  we   let
 $\gamma_b=D'_{\End E}(jj^*)$. By Lemma~\ref{lmalma} and \eqref{projj}
we then have 
\begin{equation}\label{btva}
\tilde\gamma_b=(Df-f\wedge\phi)\w\sigma
\end{equation}
and %
\begin{equation}\label{btre}
\dbar\tilde\gamma_b=
(Df-f\w\phi)\w\dbar\sigma+(\Theta f+f\w \dbar\phi)\w\sigma.
\end{equation}\label{ntre}
The following formula is the key point in the analysis of the  singularities
of $c(D_Q)$.

\begin{prop}\label{cqprop}
In   $X\setminus Z$ we have the explicit formula 
\begin{equation}\label{cq}
c(D_Q)=
\int_e f\w\sigma \w e^{\tilde I+\al\tilde\Theta
-\al Df\w\dbar\sigma}.
\end{equation}
\end{prop}


\begin{proof}
Since $\Theta_b=\Theta-\dbar\gamma_b$ we have by \eqref{btre} that
$$
\widetilde\Theta_b=\widetilde\Theta-\big((Df-f\w\phi)\w\dbar\sigma+(\Theta f+f\w \dbar\phi)\w\sigma\big).
$$
For any section $A$ of  $\End(E)$,  
\begin{equation}\label{guppi}
\int_e f\w\sigma\w \tilde A_{m-1}=\int_e f\w \sigma \w e^{\tilde A}
\end{equation}
is the determinant of the restriction of $A$ to $Q$, that is,  the determinant
of $gAg^*$.
In view of \eqref{thb}  therefore  the expression on the right hand side of 
\eqref{cq} is equal to $\det(I_Q+\al \Theta_Q)=c(D_Q)$.
\end{proof}

Now, let $v$ and $b$ be the forms in $X\setminus Z$ defined by
\eqref{gaas} and \eqref{bla}.

\begin{prop}\label{thmfyra}
\noindent(i)   The forms $v$, $b$, $c(D_Q)$,  and $c(D_S)\wedge c(D_Q)$
are  locally integrable in $X$. 

\smallskip
\noindent(ii)  If the  natural extensions are denoted by
capitals, then 
\begin{equation}\label{BV}
-2\al \partial V=B,
\end{equation}
and 
\begin{equation}\label{CCC}
-dd^cV=c(D_E)-C(D_S)C(D_Q).
\end{equation}
\end{prop}

\begin{proof}
This is clearly a local question at $Z$. Locally
we can write $f= f_1 e_j+\cdots +f_m e_m$, where
$e_j$ is a local holomorphic frame for  $E$.
In a small neighborhood $U$ of a given point
in $X$, Hironaka's theorem provides an  $n$-dimensional complex manifold $\widetilde U$
and a proper mapping $\Pi\colon\widetilde U\to U$ which is a biholomorphism
outside $\Pi^{-1}(\{f_1\cdots f_\nu=0\})$,  and such that locally on $\widetilde U$ there
are holomorphic coordinates $\tau$ such that
$\Pi^*f_j=u^j \tau_1^{\alpha-1}\cdots\tau_n^{\alpha_n}$, where $u_j$ nonvanishing;
i.e., roughly speaking $\Pi^* f_j$ are monomials.
By a  resolution over a  suitable toric manifold,
following  \cite{BY1} and  \cite{PTY}, we 
may  assume in the same way that one of the functions so obtained
divides  the other
ones. For simplicity we will make a slight abuse of notation
and suppress all occurring $\Pi^*$ and thus denote these functions by $f_j$ as well.
We  may therefore assume that $f=f_0f'$ where $f_0$ is a holomorphic
function and $f'$ is a non-vanishing section. 
Since $\sigma=f^*/|f|^2$,  it follows that    
$
\sigma=\sigma'/f_0
$
where $\sigma'$ is smooth, and hence 
$$
\widetilde{jj^*}=f\w\sigma=f'\w\sigma'
$$ 
is smooth in this resolution.
Moreover, $Df\w\dbar\sigma=Df'\w\dbar\sigma'+\cdots$, where $\cdots$ denote  terms that contain
some factor $f'$ or $\sigma'$. In view of
Proposition~\ref{cqprop}  it follows that (the pullback of) $c(D_Q)$ is smooth,
and therefore locally integrable. Since the push-forward of a locally integrable form is
locally integrable we can conclude that $c(D_Q)$ is locally integrable.

\smallskip
It follows that also 
 $\tilde\gamma_b=D'(f\w\sigma)$ and $\dbar\tilde\gamma_b$ are smooth.
Since \eqref{BV} and  \eqref{CCC} hold in $X\setminus Z$  and 
 $c(D_E)$ is smooth,  it follows that all the forms are  smooth
in the resolution. 
We can conclude that all the forms are locally integrable in $X$ and that
\eqref{BV}  and \eqref{CCC} hold.
\end{proof}

The presence of the factor $\widetilde{jj^*}=f\wedge\sigma$ implies that,
cf.,  \eqref{bla},
\begin{equation}\label{vformel}
v=\sum_{\ell=1}^{m-1}\frac{(-1)^{\ell}}{2\ell}
\int_e f\wedge\sigma\wedge (\tilde I+\al\tilde\Theta-\al Df\w\dbar\sigma)_{m-1-\ell}
\wedge (-\al Df\w\dbar\sigma)_\ell.
\end{equation}

\smallskip
\begin{df}\label{plask}
We   define the current $W$ as
\begin{multline}\label{wdef}
W=\log(1/|f|)c(D_Q)-V=\\
\log(1/|f|)\int_ef\w\sigma\w(\al\tilde\Theta+\tilde I-\al Df\w\dbar\sigma)_{m-1}\\
-\sum_{\ell=1}^{m-1}\frac{(-1)^{\ell}}{2\ell}
\int_e f\wedge\sigma\wedge (\tilde I+\al\tilde\Theta-\al Df\w\dbar\sigma)_{m-1-\ell}
\wedge (-\al Df\w\dbar\sigma)_\ell.
\end{multline}
\end{df}

In particular, if $E$ is a line bundle, i.e., $m=1$, then $V=0$, and since  $\sigma\cdot f=1$
we have that $W=\log(1/|f|)$.
It is now a simple matter to conclude the proof of Theorem~\ref{thmett}.

\begin{proof}[Proof of Theorem~\ref{thmett}]
Consider  a resolution of singularities in which  $f=f_0 f'$ with $f'$ non-vanishing,
as in the proof of Proposition~\ref{thmfyra}.
Then we know that  $c(D_Q)$ is smooth, and therefore 
$\log|f| c(D_Q)$ is locally integrable there. 
Moreover, since $\log|f|=\log|f_0|+\log|f'|$ we have that
\begin{multline*}
\lambda|f|^{2\lambda}\frac{\al\partial|f|^2\w\dbar|f|^2}{|f|^4}\w c(D_Q)=\\
\lambda|f_0|^{2\lambda}|f'|^{2\lambda}
\al\Big(\frac{df_0}{f_0}+\frac{\partial|f'|^2}{|f'|^2}\Big)\w
\Big(\frac{d\bar f_0}{\bar f_0}+\frac{\dbar|f'|^2}{|f'|^2}\Big)\w c(D_Q).
\end{multline*}
 This form is locally integrable  for $\lambda>0$ and tends to
$$
[f_0=0]\w c(D_Q)=dd^c(\log|f|c(D_Q)){\bf 1}_{\{f_0=0\}}
$$
when $\lambda\to 0$, where $[f_0=0]$ is the current of integration over the divisior defined by $f_0$.
Thus $M$ is a closed current of bidegree $(*,*)$ and order zero in $X$ with support on $Z$. 
Thus, see, e.g., \cite{Dem}, $M_k=0$ for $k<p=\codim Z$ and $M_p=\sum_j \alpha_j Z^p_j$
for some  numbers  $\alpha_j$. To see that $\alpha_j$ is  precisely the multiplicity
of $f$ on $Z_j^p$ we can locally deform the Hermitian metric to a trivial metric.
Then $\Theta=0$ and a straight-forward computation, see \cite{A3}, reveals that
$c_{p-1}(D_Q)=(dd^c\log|f|)^{p-1}$. Therefore,
$M=dd^c(\log|f|(dd^c\log|f|)^{p-1})$ which is equal to the multiplicity times
$[Z_j^p]$ according to King's formula, see \cite{King} and \cite{Dem}.
Thus part $(i)$ of the theorem is proved. 
Since  $c(D_S)-1=c_1(D_S)=dd^c\log(1/|f|)$ we have
$$
dd^c(\log(1/|f|) c(D_Q))= C(D_S)\w C(D_Q)-C(D_Q) -dd^c(\log|f| c(D_Q)){\bf 1}_Z.
$$
Now  part $(ii)$ follows from  Proposition~\ref{thmfyra}, cf, \eqref{wdef}.
\end{proof}

\section{A direct approach to \eqref{Apa}}\label{gratta}
We use the same notation as in the previous section. 
In \cite{BB},  Berndtsson introduced the deformation   $D_a=D-\gamma_a$  of $D$ on $E$, where
\begin{equation}\label{gammaa}
\tilde\gamma_a= Df\w\sigma,
\end{equation}
in order to construct Koppelman  formulas for $\dbar$ on manifolds. 
He proved formula \eqref{aek} below for $k=m$
(i.e., $\dbar a_m=d a_m=c_m(E)$).
For the general case
first we  must understand the geometric meaning of $D_a$.
Since   $D_a f=0$,  we have that  $D_a\xi$ is in $S$ if
$\xi $  is a section of $S$.
Moreover, if $\xi$ is a section of   $S^{\perp}$, then
$D_a\xi= D_E\xi$.
Now
\begin{equation}\label{konn}
g\xi\mapsto g(D_a\xi)
\end{equation}
is a well-defined connection on $Q$, and we claim that it is actually
the Chern connection $D_Q$.
In fact, if $\eta=g\xi$, then
$$
D_Q\eta=g(D_E(g^*\eta))=g(D_a(g^*\eta))=g(D_a\xi).
$$
It follows that
$\Theta_Q\eta=g(\Theta_a\xi)$,
and since $\Theta_a\xi=0$ if $\xi$ takes values in $S$, we have that
\begin{equation}\label{vin}
\al\Theta_a\sim
\left(\begin{array}{cc}
0 & * \\
0 & \al\Theta_{Q}
\end{array}
\right)
\end{equation}
with respect to the smooth isomorphism \eqref{iso}. Therefore,
$$
\al\Theta_a+I_E\sim
\left(\begin{array}{cc}
I_S & * \\
0 & I_{Q}+\al\Theta_{Q},
\end{array}
\right),
$$
and taking the determinant, we find that
\begin{equation}\label{kraft}
c(D_Q)=c(D_a).
\end{equation}

\begin{prop}\label{rost}
If $\gamma_a$ is  defined by \eqref{gammaa}, then
\begin{equation}
-tD\tilde\gamma_a+t^2\widetilde{\gamma_a\w\gamma_a}=
-t(Df\w\dbar\sigma+\Theta f\w\sigma)+
(t-t^2)Df\w\phi\w\sigma.
\end{equation}
\end{prop}

\begin{proof}
A simple computation yields 
$$
D\tilde\gamma_a=\Theta f\w\sigma+Df\w\dbar\sigma +
Df\w\phi\w\sigma
$$
and 
$$
\widetilde{\gamma_a\w\gamma_a}= Df\w\sigma\cdot Df\w\sigma,
$$
where the dot means the natural contraction of $E$ and $E^*$
so that $\xi\cdot (\alpha\w\eta)=\alpha (\xi\cdot\eta)$ if
$\xi$ and $\eta$ are sections of  $E$ and $E^*$, respectively, and $\alpha$ is a form.
Since 
$\sigma\cdot Df=-D'\sigma\cdot f=\phi$
we get the  desired formula.
\end{proof}

\begin{prop}\label{aprop}
If 
\begin{equation}\label{aformel}
a=\int_e 
\al Df\w\sigma\w e^{\tilde I+\al\tilde\Theta}\w
\sum_{\ell=0}^\infty
\frac{(-\al Df\w\dbar\sigma)^{\ell}}{(\ell+1)!}
\end{equation}
then
\begin{equation}\label{aek}
\dbar a= d a=c(D_E)-c(D_Q)
\end{equation}
in $X\setminus Z$. 
\end{prop}


\begin{proof}
We choose the homotopy  $D_t=D-t\gamma_a$ 
between $D=D_0$ and $D_1=D_a$.
In view of \eqref{smuts}, \eqref{lapp},  and Proposition~\ref{rost}  we  have that
$$
a=\int_e\int_0^1 \al Df\w\sigma\w e^{\tilde I+\al\tilde\Theta-
t\al (\Theta f\w\sigma+Df\w\dbar\sigma) -(t-t^2)Df\w\phi\w\sigma}dt
$$
satisfies the second equality in \eqref{aek} in $X\setminus Z$. 
Noticing that $\sigma\w\sigma=0$, a computation of the $t$-integral yields
\eqref{aformel}.
Since $a$ has bidegree $(*,*-1)$ and $da$ has bidegree $(*,*)$
it follows that $\dbar a=da$.
\end{proof}

The forms $a$ and $b$ are related in the following way.

\begin{prop}\label{abas}
In $X\setminus Z$ we have that 
\begin{equation}\label{aba}
b=a+\al\partial\log|f|^2\w c(D_Q)
\end{equation}
\end{prop}

\begin{proof}
Starting with  \eqref{gaas} we have 
\begin{multline*} 
b=\int_e\al(Df-f\w\phi)\w\sigma\w e^{\tilde I+\al\tilde\Theta}
\w\sum_{\ell= 0}^\infty\frac{(-\al Df +\al f\w\phi)^\ell}{(1+\ell)!}\w(\dbar\sigma)^\ell=\\
-\int_e e^{\tilde I+\al\tilde\Theta}\w\sum_{\ell=0}^\infty\frac{(-\al Df+\al f\w\phi)^{\ell+1}}
{(\ell+1)!}\wedge \sigma\w(\dbar\sigma)^\ell=\\
-\int_e
e^{\tilde I+\al\tilde\Theta -\al Df+\al f\wedge\phi}\w \sum_{\ell=0}^\infty
\sigma\w(\dbar\sigma)^\ell = \\
-\int_e e^{\tilde I+\al\tilde\Theta-\al Df}
\wedge(1 +\al f\w\phi)\wedge  \sum_\ell\sigma\w(\dbar\sigma)^\ell.
\end{multline*}
In view of  \eqref{batt} and \eqref{cformel},  recalling that
  $\phi=-\partial\log|f|^2$, we now get \eqref{aba}.
\end{proof}

By a resolution of singularities as in the proof of Proposition~\ref{thmfyra}  above one can 
see that $a$ is locally integrable. Let $A$ denote its natural extension. 
By such a resolution  one 
can also verify that the formal computation (using Proposition~\ref{abas})
$-2\al\partial (\log(1/|f|)c(D_Q)-V)=B-\al\partial\log|f|^2\w C(D_Q)=A$
is ligitimate, and thus we have
\begin{equation}\label{skott2}
A=-2\al\partial W.
\end{equation}
As a consequence we get that $\dbar A=dA=c(D_E)-c(D_Q)-M$.

\section{Factorization of currents}\label{faktor}

Since $a$ and $c(D_Q)$ are  locally integrable, 
 $|f|^{2\lambda} a$  and $|f|^{2\lambda}c(D_Q)$ are well-defined currents  for 
 $\Re\lambda>-\epsilon$ and we have 
\begin{equation}\label{Af}
A=|f|^{2\lambda}a|_{\lambda=0} \quad {\rm and}\quad C(D_Q)=|f|^{2\lambda}c(D_Q)|_{\lambda=0}.
\end{equation}
It also follows that
\begin{equation}\label{Mf}
M=-d|f|^{2\lambda}\w a|_{\lambda=0}=-\dbar|f|^{2\lambda}\w a|_{\lambda=0}.
\end{equation}
Now consider the expression \eqref{aformel} for $a$.
Since each term in 
$\exp(\tilde I+\al\tilde\Theta)$ has the same degree in
$e_j$ and $e^*_k$ it  must be  multiplied by
terms with the same property in order to get a product  with full degree.
Therefore we can rewrite $a$ as 
\begin{equation}\label{batt}
a=-\int_e e^{\tilde I+\al\tilde\Theta-\al Df}\w
\sum_0^\infty \sigma\w (\dbar\sigma)^{\ell}.
\end{equation}
In  \cite{A2} we introduced the currents
$$
U=|f|^{2\lambda}\frac{\sigma}{1-\dbar\sigma}\big|_{\lambda=0}=
|f|^{2\lambda}\wedge\sigma\wedge \sum_\ell 
(\dbar\sigma)^{\ell-1}\big|_{\lambda=0}
$$
and
$$
R=\dbar|f|^{2\lambda}\wedge\frac{\sigma}{1-\dbar\sigma}\big|_{\lambda=0}=
\dbar|f|^{2\lambda}\wedge\sigma\wedge \sum_\ell (\dbar\sigma)^{\ell-1}
\big|_{\lambda=0}.
$$
It is part of the statement that the right hand sides are current valued holomorphic
functions for $\lambda>-\epsilon$,   evaluated at $\lambda=0$. 
In general  $U$ and $R$ are {\it not} locally integrable.
The current  $R$ is supported on $Z$,
$$ 
R=R_p+\cdots +R_{\min(m,n)},
$$ 
where $R_k$ is the component of bidegree $(0,k)$ taking values in $\Lambda^k E^*$,
 and $(\delta_f-\dbar)U=1-R $.
In view of \eqref{batt}, \eqref{Af}, and \eqref{Mf}  we have the factorization
formulas 
\begin{equation}\label{gris}
M=\int_e  e^{\al\widetilde\Theta+\tilde I-\al Df}\w R,
\end{equation}
\begin{equation}\label{gris2}
A=-\int_e  e^{\al\widetilde\Theta+\tilde I-\al Df}\w U,
\end{equation}
and moreover, cf. \eqref{cq}, 
\begin{equation}\label{cformel}
C(D_Q)= \int_ef\wedge\sigma\wedge e^{\al\tilde\Theta+\tilde I-\al Df\wedge\dbar\sigma}
=\int_e   e^{\al \widetilde\Theta+\tilde I-\al Df}\w f\w U.
\end{equation}

\section{Positivity}\label{poss}

Let $E\to X$ be a Hermitian  holomorphic bundle as before and let $e_j$ be an orthonormal
local frame. A  section
$$
A=i\sum_{jk} A_{jk}\otimes e_j\otimes e_k^*
$$
of   $T^*_{1,1}(X)\otimes\End(E)$ is Hermitian  if $A_{jk}=-\overline{A_{kj}}$.
It then induces a Hermitian  form $a$ on $T^{1,0}(X)\otimes E^*$ by
$$
a(\xi\otimes e_j^*,\eta\otimes e_k^*)= A_{jk}(\xi,\bar\eta),
$$ 
if  $\xi,\eta$ are $(1,0)$-vectors. We say that $A$ is (Bott-Chern) positive, $A\ge_B 0$ if
the form $a$ is positively semi-definite. 
In the same way any Hermitian  $A$ induces a Hermitian  form $a'$ on $T^{1,0}(X)\otimes E$
and it is called Nakano positive, $A\ge_N 0$,  if $a'$ is positively semi-definite.

Notice that 
$\al\Theta$ is Hermitian;  it is said to be Nakano positive if 
$\al\Theta\ge_N 0$. Analogously we say that $E$ is positive,
$E\ge_B 0$,  if $\al\Theta \ge_B 0$. 
Neither of these positivity concepts implies the other one unless $m=1$.

Since $\Theta_{jk}(E^*)=-\Theta_{jk}(E)$ it follows that
$E$  is positive in our sense if and only if $E^*$ is Nakano negative.
The next proposition explains the interest of Bott-Chern positivity
in this context.

\begin{prop}
Let 
\begin{equation}\label{seq}
0\to S\to E\to  Q\to 0
\end{equation}
 be a short exact sequence of Hermitian 
holomorphic vector bundles. Then
$E\ge_B 0$ implies that $Q\ge_B 0$.
\end{prop}

\begin{proof}
It is well-known, see for instance \cite{Dem},  that $E\le_N 0$ implies
that $S\le_N 0$. From  the  sequence \eqref{seq}  above we get  the  exact sequence
$0\to Q^*\to E^*\to S^*\to 0$. 
Since $E^*\le_N 0$ implies $ Q^*\le_N 0$,  it follows that 
$E\ge_B 0$ implies $Q\ge_B 0$.
\end{proof}

The next simple lemma reveals that our definition of Bott-Chern positivity  coincides
with the one used in \cite{BC}.

\begin{lma}
$A\ge_B 0$ if and only if there are sections $f_\ell$ of $T^*_{1,0}(X)\otimes E$ such that
\begin{equation}\label{apa}
A=i\sum_\ell f_\ell\otimes f^*_\ell.
\end{equation}
\end{lma}

Observe  that if  $f_\ell=\sum f_j^\ell\otimes e_j$, then $f^*_\ell=\sum \bar f^\ell_j\otimes e_j^*$ 
since $e_j$ is ortonormal.

\begin{proof}
If \eqref{apa} holds, then 
$$
a(\xi,\xi)=\sum_\ell f_\ell(\xi)f_\ell^*(\xi^*)=\sum |f_\ell(\xi)|^2\ge 0
$$
for all $\xi$ in $T^{1,0}\otimes E^*$.
Conversely, if $a$ is positive, it is diagonalizable, and so there is a basis
$f_\ell$ for $T^*_{1,0}\otimes E$ such that \eqref{apa} holds.
\end{proof}

If we identify $f_\ell$ with $\sum f^\ell_j\w e_j$ as before, then \eqref{apa} means that 
\begin{equation}\label{amojs}
\tilde A= -i\sum_\ell f_\ell\w f^*_\ell.
\end{equation}

\smallskip

If $B=\sum B_{jk} e_j\otimes e_j^*$ is a scalar-valued section of  $\End E$, then 
it is Hermitian  if and only if $B_{jk}=\bar B_{kj}$ and it is positively
semi-definite if and only if
$$
B=\sum_\ell g_\ell\otimes g^*_\ell
$$
for some sections $g_\ell$ of  $E$; or equivalently, 
\begin{equation}\label{bmojs}
\tilde B=\sum_\ell g_\ell\wedge g^*_\ell.
\end{equation}

\begin{prop}\label{skuta}
Assume that  $A_j$ are $(1,1)$-form-valued Hermitian  sections of  $E$ and
$B_k$ scalarvalued sections, such that $A_j\ge_B 0$ and $B_k\ge 0$. Then
\begin{equation}\label{gott}
\int_e\tilde A_1\w\ldots\w\tilde A_r\w\tilde B_{r+1}\w\ldots\w\tilde B_m
\end{equation}
is a positive $(r,r)$-form.
\end{prop}

\begin{proof}
In view of \eqref{amojs} and \eqref{bmojs}, we see that \eqref{gott} is a sum of
terms like
\begin{multline*}
\int_e(-i)^r f_1\w f_1^*\w\ldots\w f_r\w f_r^*\w g_{r+1}\w g^*_{r+1}\w\ldots
\w g_m\w g^*_m = \\
(-i)^rc_{m-r}\int_e f_1\w\ldots f_r\w\ldots g_m\w f^*_1\w\ldots\w f^*_r\w\ldots g^*_m=\\
(-i)^r c_{m-r}\int_e \omega\w e_1\w\ldots\w e_m\w\bar\omega\w e_1^*\w\ldots\w e^*_m,
\end{multline*}
where $\omega$ is an $(r,0)$-form and $c_p=(-1)^{p(p-1)/2}=i^{p(p-1)}$.
By further simple computations, 
\begin{multline*}
(-i)^rc_{m-r}(-1)^{mr}\int_e\omega\w\bar\omega\w e_1\w\ldots\w e_m\w e_1^*\w\ldots\w e^*_m=\\
(-i)^rc_{m-r}(-1)^{mr}c_m\omega\w\bar\omega=i^{r^2}\omega\w\bar\omega
\end{multline*}
the proposition follows, since the last form is positive.
\end{proof}

\begin{prop}
If $E\ge_B 0$ (or $E\ge_N 0$), then the Chern forms
$c_k(D_E)$ are positive for  all $k$.
\end{prop}

\begin{proof}
Since $\alpha\Theta\ge_B 0$ by assumption, and clearly $I\ge 0$, it follows from
Proposition~\ref{skuta} that 
$$
c_k(D_E)=\int_e (\al\widetilde\Theta)_k\w\tilde I_{m-k}
$$
is positive.
\end{proof}

\begin{proof}[Proof of Theorem~\ref{thmettx}]
We have just seen that $c(D_E)\ge 0$.
From \eqref{Md} it follows   that   the current
$M_k$ is positive if $c_{k-1}(D_Q)$ is positive.
From  \eqref{cq} we have that 
\begin{multline}\label{ckform}
c_{k-1}(D_Q)= \int_e f\w\sigma\w (\al\widetilde\Theta-\al Df\w\dbar \sigma)_{k-1}\w\tilde I_{m-k}=\\
\sum_{j=1}^{k-1}\int_e f\w\sigma\w (\al\widetilde\Theta)_{k-1-j}\wedge
(-\al Df\w\dbar \sigma)_j\w\tilde I_{m-k}.
\end{multline}
If $s=f^*$ as before, then  $\sigma=s/|f|^2$, and 
therefore we have  
\begin{equation}\label{ormskinn}
c_{k-1}(D_Q)=\sum_{j=1}^{k-1}\int_e \frac{f\w s}{|f|^2}\w 
\Big(\frac{-\al Df\w\dbar s}{|f|^2}\Big)_j\w
(\al\widetilde\Theta)_{k-1-j}\wedge\tilde I_{m-k}.
\end{equation}
Since $\dbar s=(Df)^*$ it now follows immediately from Proposition~\ref{skuta}
that $c_k(D_Q)$ is positive if $\al\Theta\ge_B 0$.

\smallskip
It remains to see that one can choose $W$ so that it is positive where
$|f|<1$. 
Notice that if some of the $A_j$ in \eqref{gott} are replaced by
$A'_j\ge_B A_j$, then the resulting form will be larger;
this follows immediately from the proof. 
Now, 
$\log(1/|f|)c(D_Q)$ is positive when  $|f|<1$.
From \eqref{vformel}  we have that
$$
v_k=\sum_{\ell=1}^k \frac{(-1)^\ell}{2\ell}\int_e f\w\sigma\w(\al\widetilde\Theta-\al Df\w\dbar\sigma)_{k-\ell}
\w(-\al Df\w\dbar\sigma)_\ell\w\tilde I_{m-k-1}.
$$
Since this  is an alternating sum of positive terms  it has no 
sign. If we replace each factor  $ -\al Df\w\dbar\sigma$
by $\al\widetilde\Theta-\al Df\w\dbar\sigma$, then we get a larger
form which in addition is closed, since it is just
a certain constant times $c_k(D_Q)$, cf., \eqref{ckform}.
Therefore, for a suitable constant $\nu_k$
$-v'_k=-v_k+\nu_k c_k(D_Q)$ is  a positive form and
$dv_k'=dv_k$.  Thus the current
$$
W'_k= -V_k+\nu_kC_k(D_Q)+\log(1/|f|)C_k(D_Q)
$$
will have the stated property.
\end{proof}

The modification of $v$ in last part of the proof is precisely as in \cite{BC}
but with  our notation,  and for an
arbitrary $k$ rather than just $k=m-1$.
It is not necessary to consider each $v_k$ separately.
By the same argument one can see directly that
$-v'=-v+\nu c(D_Q)$ is positive   if $\nu$ is appropriately chosen,
and $dv'=dv$.

\smallskip

One can  prove that
if we multiply \eqref{ormskinn} with $\lambda \partial|f|^2\w\dbar|f|^2/|f|^2$
and let $\lambda\to 0^+$, then all terms  with $j<p-1$  will disappear;
see for instance the proof of Theorem~1.1 in \cite{A2}.
We thus have

\begin{prop}\label{mprecis}
If $p=\codim\{f=0\}$, then
\begin{multline*}
M_k=\lim_{\lambda\to 0^+}\lambda|f|^{2\lambda}\al\frac{\partial|f|^2\w\dbar|f|^2}{|f|^2}
\wedge \\
\sum_{j=p-1}^{k-1}\int_e \frac{f\w s}{|f|^2}\w 
\Big(\frac{-\al Df\w\dbar s}{|f|^2}\Big)_j\w
(\al\widetilde\Theta)_{k-1-j}\wedge\tilde I_{m-k}.
\end{multline*}
\end{prop}

From this formula it is apparant that
 $M_k$ vanishes if $k<p$, and that $M_p$ is positive, regardless
of $\al\Theta$. One can also derive  this formula from
\eqref{gris}.

\begin{remark}
When $k>p$, $M_k$ depends on the metric, but there is still
a certain uniqueness:
Let  $Z^k$ be   the union of the irreducible components $Z_j^k$ of $Z$ of
codimension $k$. One can verify,  see \cite{A3}, that
the restriction of $M_k$ to $Z^k$ is a sum
$$
\sum_j \alpha_j^k[Z_j^k],
$$
where  $\alpha_j^k$ are nonnegative numbers that are independent of 
the metric. However the geometric meaning of these numbers
is not clear to us.
\end{remark}

\section{Some examples}

The first two examples suggest  that not only the component $M_p$ of the current
 $M$ is of interest.

\begin{ex}\label{snitt}
Let us  assume that $X$ is compact, and that 
we  have sections $f_j$ of  rank $m_j$ bundles $E_j\to X$,
such that $\sum m_j=n$. If $E=\oplus E_j$ and $f=(f_1,\ldots, f_r)$, then 
the intersection number $\nu$  of the varieties $Z_j=\{f_j=0\}$ is equal to the integral
of 
$$
c_n(E)=c_{m_1}(E_1)\w\ldots\w c_{m_r}(E_r)
$$ 
over $X$. Since $M_n$ represents the cohomology class $c_n(E)$, we thus get the representation 
$$
\nu =\int_X M_n,
$$
i.e., an integral over the set-theoretic  intersection $Z=\cap Z_j$. 
If $E$ is positive then $M_n$ is positive. If $Z$ is discrete, i.e., $f$ is a complete
intersection,  then 
$M_n=[Z]$,  and in this case thus we just get 
the sum of the points in $Z$ counted with multiplicities, as expected.
\end{ex}

\begin{ex}
Let $X$ be a compact K\"ahler manifold with metric form $\omega$, and let
$f$ be a holomorphic section of  $E\to X$. If moreover $E\ge_B 0$, then 
we know that $c(D_E)$, $M$, and $c(D_Q)0$ are all positive. 
Because of \eqref{bas}, we  therefore have that
$$
\int_X M_k\w\omega_{n-k}=\int_X c_k(D_E)\w\omega_{n-k}-\int_X c_k(D_Q)\w\omega_{n-k}\le
\int_X c_k(D_E)\w\omega_{n-k}.
$$
Thus we get an upper bound of the total mass of $M_k$ in terms of the Chern class $c_k(E)$.
Taking $k=p=\codim Z$ we get the estimate
$$
{\rm area}(Z^p)=\int_X[Z^p]\le \int_X c_p(E)\w\omega_{n-p}.
$$
\end{ex}

\begin{ex}\label{om}
Now assume that  $X=\P^n$, let 
$$
\omega=\al \partial\dbar\log|z|^2=dd^c\log|z|
$$
denote the Fubini-Study metric and notice that
$$
\int_{P^n}\omega^n=1,
$$
that is, the total area of $\P^n$ is $1/n!$.

\smallskip
Assume that $F_1,\ldots, F_m$ are polynomials in $\C^n$ which
form a complete intersection. If $F_j$ has degree $d_j$ (depending on  $z'=(z'_1,\ldots,z_n')$) then
the  the homogenization  $f_j(z)=z_0^{d_j}F(z'/z_0)$ is a 
$d_j$-homogeneous polynomial in $\C^{n+1}$ and hence corresponds to a section
of the line bundle $\O(d_j)\to\P^n$. Thus $f=(f_1,\ldots,f_m)$ is a section
of $E=\oplus \O(d_j)$. 
If $E$ is equipped with the natural metric, i.e., 
$$
\|h([z])\|^2=\sum_j \frac{|h(z)|^2}{|z|^{2d_j}}
$$
for a section $h=\oplus h_j$ of $E$
(here $[z]$ denotes the point on $\P^n$ corresponding to the point $z\in \C^{n+1}\setminus\{0\}$
under the usual projection),
then
it  is easy to check that $E\ge_B 0$.
Therefore $M_m\ge 0$, and since moreover,
$$
M_m|_{\C^n}=[Z],
$$
if $Z$ here denotes the zero variety $\{F=0\}$ in $\C^n$, then
$$
{\rm  area (Z)}=\int_{C^n}[Z]\w\omega_{n-m}\le
\int_{\P^n}M_m\w\omega_{n-m}=
\int_{\P^n}{c_m(D_E)}\w\omega_{n-m},
$$
since $c_m(D_Q)=0$. Here ``area''  refers to the projective area of course.
 However,
$
c(D_E)=(1+d_1\omega)\w\ldots\w(1+d_m\omega),
$
and so
$$
c_m(D_E)=d_1\cdots d_m \omega^m.
$$
Hence 
$$
{\rm  area (Z)}\le d_1\cdots d_m \frac{1}{(n-m)!}.
$$
We also notice that the deviation from equality is precisely the
total mass of $M_m$ on the hyperplane at infinity.
If   $m=n$  we get Bezout's theorem 
$$
\#\{F=0\}\le d_1\cdots d_n.
$$
\end{ex}


\begin{ex}
If  $f$ is a complete intersection, i.e., $p=m$, and $W_{m-1}$ denotes the
component of bidegree $(m-1,m-1)$, then 
$$ %
dd^cW_{m-1}=c_m(D_E)-[Z];
$$ %
this  means that $W_{m-1}$ is  a Green current for the cycle $Z=\sum \alpha_j Z_j$.

In the case when $E=L_1\oplus\cdots\oplus L_m$
for some line bundles $L_k$,
hence $c_m(D_E)=c_1(D_{L_1})\w\ldots\w c_m(D_{L_m})$,
and $f=(f_1,\ldots,f_m)$, where $f_j$ are holomorphic sections of $L_j$,
such a Green current was obtained already in   \cite{BY1}.
\end{ex}

\begin{ex}\label{koppelman}
Let $X$ be  a compact manifold such that there is a holomorphic section $\eta$ of  some
vector bundle $H \to X\times X$ that defines the diagonal $\Delta\subset
X\times X$; for instance $X$ can be complex projective space. 
From  Theorem~\ref{thmett} we get a current $W_n$ such that
$dd^c W=c_n(D_H)-[\Delta]$. If we let $K(\zeta,z)=-W_n$ and
 $P(\zeta,z)=c_n(D_H)$, then 
$$
dd^cK=[\Delta]-P,
$$
and this leads to the   Koppelman type formula
\begin{multline}\label{kman}
\phi(z)-\int P(\zeta,z)\w\phi(\zeta)=\\
dd^c\int_X K\w \phi -d\int_X K\w d^c\phi
+d^c\int_X  K \w d\phi+
\int_X K\w dd^c\phi
\end{multline}
for the $dd^c$-operator.
In particular, if $\phi$ is closed $(k,k)$-form such that $d\phi=0$, then  $d^c\phi=0$ as well, 
and thus 
$$
v=\int_X K\w\phi
$$
is  an explicit solution to
$dd^c v=\phi-\int P\w\phi$. 
However if $X$ is non-compact one gets  boundary integrals.
It would be desirable  to refine the construction to include somehow 
an appropriate line bundle with a metric that vanishes at the boundary,
in order  to obtain $dd^c$-formulas for, say,  domains in $\C^n$.
\end{ex}

\begin{ex}\label{greenex}
Assume that  $f$ is a holomorphic section of some Hermitian bundle
$E\to X$ with zero variety $Z$. If $f$ is locally a complete intersection
we have  seen that the current $W_{m-1}$ from Theorem~\ref{thmett}
is a Green current for $[Z]$. In general we  have that
$dd^cW_{p-1}=c_p(D_E)-c_p(D_Q)-[Z^p]$ so we only get a current $w$ such that
$dd^cw=[Z^p]-\gamma$, where $\gamma$ is locally integrable. 
However, there is another and simpler way  to find such a current  $w$,  due to
Meo, \cite{Meo2}. 
\begin{prop}[Meo]
Let $f$ be a holomorphic section of a Hermitian vector bundle $E\to X$.
The forms  
$$
w=\log|f|\big (dd^c\log|f|)^{p-1}{\bf 1}_{X\setminus Z}\big)
$$
and 
$$
\gamma=-(dd^c\log|f|)^p{\bf 1}_{X\setminus Z}
$$
are locally integrable on $X$ and 
\begin{equation}\label{melog}
dd^cw=[Z^p]-\gamma. 
\end{equation}
\end{prop}

For the reader's convenience we provide a simple proof
based on Hironaka's theorem.

\begin{proof}[Sketch of proof]
Let  $f=f_0 f'$ be as before, i.e.,  $f_0$ is holomorphic and $f'$ is a non-vanishing section. Then
$\log|f|=\log|f_0|+\log|f'|$,  and hence $dd^c\log|f'|$ is smooth and
$dd^c\log|f_0|=[f_0=0]$ has support on  the inverse
image $\tilde Z$ of $Z$ in the resolution. Thus 
$$
w=(\log|f_0|+\log|f'|)(dd^c\log|f'|)^{p-1}, \quad \gamma=(dd^c\log|f'|)^p
$$
are both locally integrable in the resolution and hence also on the original manifold.
Moreover,
$$
dd^c w=[f_0=0]\wedge(dd^c\log|f'|)^{p-1}+\gamma,
$$
in particular $(dd^c w){\bf 1}_{\tilde Z}$ is closed, and hence
$T=(dd^c w){\bf 1}_Z$ is a closed current on $X$ of order zero.  
Therefore $T=\sum\alpha_j[Z_j^p]$,  and since we can deforme the metric into a trivial metric
locally,   it follows from King's formula, \cite{King},
that $\alpha_j$ are precisely the multiplicities of $f$ on $Z_j^p$.
\end{proof}

Assume now that $X$ is a compact manifold such that there exists a holomorphic section
$\eta$ of some Hermitian bundle $H\to X\times X$ as in Example~\ref{koppelman} above. 
If furthermore the kernel $K$ is reasonably regular  we can assume that
$$
dd^c\int_\zeta K(\zeta,z)\w \psi(\zeta)=\psi(z)-\int_\zeta P(\zeta,z)\w\psi(\zeta)
$$
for any $(k,k)$-current $\psi$. We then have the explicit solution 
$$
g= w+\int_\zeta K(\zeta,z)\w\gamma(\zeta)
$$
to the  Green equation
$
dd^c g= [Z^p]-\alpha,
$
where $\alpha$ is the smooth form
$$
\alpha=\int_\zeta P(\zeta,z)\w\gamma(\zeta).
$$
\end{ex}

The last example will be elaborated in a forthcoming paper.

\def\listing#1#2#3{{\sc #1}:\ {\it #2},\ #3.}

\end{document}